\documentclass{amsart}
\usepackage{amssymb}
\usepackage[T1]{fontenc}
\usepackage{enumerate}
\usepackage{url}

\newcommand\N{\mathbb N}
\newcommand\Z{\mathbb Z}
\newcommand\Q{\mathbb Q}
\newcommand\R{\mathbb R}

\newcommand\ph\varphi
\newcommand\ps\psi
\newcommand\ep\varepsilon
\newcommand\rh\varrho
\newcommand\al\alpha
\newcommand\be\beta
\newcommand\ga\gamma
\newcommand\om\omega
\newcommand\ta\tau
\renewcommand\th\vartheta
\newcommand\de\delta
\newcommand\ze\zeta
\newcommand\ch\chi
\newcommand\et\eta
\newcommand\io\iota
\newcommand\la\lambda
\newcommand\si\sigma

\newcommand\Ga\Gamma
\newcommand\De\Delta
\newcommand\Th\Theta
\newcommand\La\Lambda
\newcommand\Si\Sigma
\newcommand\Ph\Phi
\newcommand\Ps\Psi
\newcommand\Om\Omega

\newtheorem{theorem}{Theorem}
\newtheorem{lemma}[theorem]{Lemma}
\newtheorem{proposition}[theorem]{Proposition}
\newtheorem{corollary}[theorem]{Corollary}
\theoremstyle{definition}
\newtheorem{definition}[theorem]{Definition}

\theoremstyle{remark}
\newtheorem{remark}[theorem]{Remark}
\newtheorem{example}[theorem]{Example}

\newcommand\x{{\bar X}}
\newcommand\zx{{\Z[\x]}}
\newcommand\rx{{\R[\x]}}
\newcommand\nnx{{\R[\x]^+}}
\newcommand\nn{{[0,\infty)}}
\newcommand\pos{{(0,\infty)}}
\renewcommand{\O}{{\mathcal O}}

\newcommand{\m}{{\mathfrak m}}
\renewcommand{\c}{{\mathcal C}}

\DeclareMathOperator{\Log}{Log}
\DeclareMathOperator{\New}{New}
\DeclareMathOperator{\init}{in}
\DeclareMathOperator{\im}{im}
\DeclareMathOperator{\Sper}{Sper}

\begin{document}
\title[Membership in archimedean semirings]
{A criterion for membership\\in archimedean semirings}
\author{Markus Schweighofer}
\address{Universit\"at Konstanz\\
         Fachbereich Mathematik und Statistik\\
         78457 Konstanz\\
         Allemagne}
\email{Markus.Schweighofer@uni-konstanz.de}
\thanks{Partially supported by the DFG project 214371
``Darstellung positiver Polynome''. The author thanks Alex Prestel for
helpful discussions}
\date{\today}
\begin{abstract}
Let $A$ be a commutative ring and $T\subset A$ a weakly divisible archimedean
semiring,
i.e., $0\in T$, $T+T\subset T$, $TT\subset T$, $\Z+T=A$ and $\frac 1r\in T$
for some integer $r\ge 2$. The classical
Real Representation Theorem says the following: If $a\in A$ satisfies
$\ph(a)>0$ for all ring homomorphisms $\ph:A\to\R$ with $\ph(T)\subset\nn$,
then $a\in T$.

The main drawback of this criterion for membership is that it is only
sufficient but far from being necessary since $\ph(a)>0$ cannot be replaced by
$\ph(a)\ge 0$ without any further conditions.
Initiated by work of Scheiderer, a lot of progress has previously been made in
overcoming this drawback but only in the case where $T$ is a preorder,
i.e., contains all squares of $A$.

A different approach enables us to prove a suitable extension of the Real
Representation theorem for the general case. If (and, to the best of our
knowledge, only if) $T$ is a
preorder, our result can easily be derived by earlier work of Scheiderer,
Kuhlmann, Marshall and Schwartz. In contrast to this earlier work, our proof
does not use and therefore shows the classical theorem.

We illustrate the usefulness of our result by deriving a theorem of Handelman
from it saying inter alia the following: If an odd power of a real
polynomial in several variables has only nonnegative coefficients, then so
do all sufficiently high powers.
\end{abstract}
\maketitle

\section{Archimedean semirings}\label{semiringsection}

Throughout this article, $A$ denotes a commutative ring. The case where the
unique ring homomorphism $\Z\to A$ (all rings have a unity and all ring
homomorphisms preserve unities) is not an embedding is formally admitted
but our results will be trivial in this case. So the reader might assume
that $A$ contains $\Z$ as a subring. Whenever we postulate that $\frac 1r$
lies in $A$ for some integer $r\ge 2$, we implicitly require that $r$
(that is $r\cdot 1$) is a unit of $A$ (i.e., invertible in $A$).

\begin{definition}\label{semiringdef}
A set $T\subset A$ is called a \emph{semiring} of $A$ if $0,1\in T$
and $T$ is closed under addition and multiplication, i.e.,
$T+T\subset T$ and $TT\subset T$. A semiring $T$ of $A$ is called a
\emph{preorder} of $A$ if it contains all the squares of $A$, i.e.,
$A^2\subset T$. We call a semiring $T$ \emph{archimedean} (with respect
to $A$) if $\Z+T=A$. We call a semiring \emph{weakly divisible} if there
is some integer $r\ge 2$ with $\frac 1r\in T$.
\end{definition}

Semirings in our sense (i.e., as subsets of rings) are often called
\emph{preprimes} (cf. \cite[Definition 5.4.1]{pd}). This goes back to
Harrison who called these objects \emph{infinite} preprimes (opposing
them to his \emph{finite} preprimes) which makes sense in a certain number
theoretic context \cite{har}. However, without the adjective ``infinite''
and in a general context, this terminology is hermeneutic. We use the term
``semiring'' and hope that other authors will follow. 
Some authors require that a semiring (``preprime'' in their terminology)
does not contain $-1$. This leads to similar technical problems than
it would to require all ideals of a ring to be proper.

This work provides a new criterion to prove membership in an archimedean
semiring. We shortly explain the well-known basic ideas. Suppose $S$ is
a compact space. In the ring $\c(S,\R)$ of continuous functions
on $S$, the nonnegative functions form an archimedean semiring $\c(S,\nn)$
(which is even a preordering).
This semiring is defined by a clear geometric property (namely being
nonnegative on the space $S$). More interesting semirings $T$, however, are
often defined in an algebraic way, for example by a set of generators. The
question arises if one can nevertheless develop criteria for membership in
$T$ of (as far as possible) geometric nature. The first step is to view
the elements of $A$ as real-valued continuous functions on some topological
space $S(T)$ (naturally associated to $A$ and $T$) such that the elements
of $T$ are nonnegative on $S(T)$.

For any semiring $T\subset A$, we set
$$S(T):=\{\ph\mid\text{$\ph:A\to\R$ ring homomorphism},
\ph(T)\subset[0,\infty)\}\subset\R^A$$
where the topology on $S(T)$ is induced by the product topology on $\R^A$,
i.e., is the weakest topology making $S(T)\to\R:\ph\mapsto\ph(a)$ for all
$a\in A$ continuous. If $T$ is archimedean, then $S(T)$ is compact
(meaning quasi-compact and Hausdorff). This is clear because
$S(T)$ equals the intersection (\ref{iplus})--(\ref{it}) appearing in the
proof of Theorem \ref{aboutissement} below. We now have a ring homomorphism
$$A\to\c(S(T),\R):a\mapsto(\ph\mapsto\ph(a))$$
sending all $a\in T$ to a function nonnegative on the whole of $S(T)$. When we
write $a$, we will often mean the image under this map. In this sense,
$\ph(a)=a(x)$ for all $x:=\ph\in S(T)$.

Often, $S(T)$ takes on a very concrete form. Concerning our motivating
example of the archimedean semiring $\c(S,\nn)\subset\c(S,\R)$ where $S$
is a compact space, it follows from basic set topology that the canonical map
$S\to S(\c(S,\nn)):x\mapsto(f\mapsto f(x))$ is a homeomorphism
(cf. \cite[4.9(a)]{gj}) allowing us to write $S=S(T)$. This illustrates the
naturality of the definition of $S(T)$. However, for this semiring $\c(S,\nn)$
our membership criterion will be inferior to the self-evident one.

Our criterion will rather be interesting in the realm of polynomials.
Throughout this article, we will consider polynomials in $n$ variables
$\x:=(X_1,\dots,X_n)$.
The polynomial ring in these $n$ variables over a commutative ring $R$ will be
denoted by $R[\x]$. For any set $P\subset\rx$, we define
$$V(P):=\{x\in\R^n\mid p(x)=0\text{\ for all $p\in P$}\}\subset\R^n.$$
Suppose that $A$ is finitely generated over a subring $R$.
Then (up to isomorphism) $A=R[\x]/I$ for some
number $n$ of variables and an ideal $I$ of $A$.
If $R\subset\R$ and $[0,\infty)\cap R\subset T$, then
every $\ph\in S(T)$ is the identity on $R$ and it is easy to see that
\begin{equation}\label{concret}
S(T)=\{x\in V(I)\mid\text{$t(x)\ge 0$ for all $t\in T$}\}\subset\R^n
\end{equation}
via the homeomorphism
$$\ph\mapsto(\ph(X_1+I),\dots,\ph(X_n+I)).$$
In particular, if $A=\rx$ and $T\subset A$ is finitely generated over
$[0,\infty)$,
say
$$T=\langle[0,\infty)\cup\{t_1,\dots,t_m\}\rangle$$
(we always write angular brackets for the generated semiring), then
$$S(T)=\{x\in\R^n\mid t_1(x)\ge 0,\dots,t_m(x)\ge 0\}\subset\R^n$$
is a so called basic closed semialgebraic set (cf. \cite[Theorem 2.4.1]{pd}).
If in addition all $t_i$ are linear (i.e., of degree $\le 1$), then $S(T)$
is a polyhedron. If this polyhedron $S(T)$ is compact (i.e., a polytope), then
it follows from a well-known theorem on linear inequalities (cf.
\cite[Theorem 5.4.5]{pd}\cite{h4}) and Proposition \ref{archprop} below
that $T$ is archimedean (and therefore also
any semiring $T'\subset A$ containing $T$). As already mentioned, the converse
is true in general: If $T$ is archimedean, then $S(T)$ is compact.

Without the linearity assumption on the $t_i\in\rx$, Schm\"udgen
\cite{sch} showed that compactness of $S(T)$ implies (and therefore is
equivalent) to the condition that the 
\emph{preordering} $T'\supset T$ generated by $t_1,\dots,t_m$, i.e.,
the semiring
$$T':=\langle\rx^2\cup\{t_1,\dots,t_m\}\rangle\supset T$$ is archimedean.

Our criterion will extend the classical criterion which is Corollary
\ref{classique} in this work. It is going back to Krivine, Stone, Kadison,
Dubois and Becker. It used to be called Kadison-Dubois theorem but due to its
(to some extent only recently revealed) complex history (see
\cite[Section 5.6]{pd}) it is now often called Real Representation Theorem.
It simply says that for a weakly divisible archimedean semiring $T\subset A$,
all $a\in A$ with $a>0$ on $S(T)$ lie in $T$.

Using what we said above, this implies for example Handelman's theorem
that any polynomial positive on a polytope is a nonnegative linear combination
of products of
the linear polynomials defining the polytope \cite[Theorem 5.4.6]{pd}\cite{h4}.
Also, it implies the corresponding weaker representation of polynomials
positive on compact basic closed semialgebraic sets proved by Schm\"udgen
\cite[Theorem 5.2.9]{pd}\cite{sch}.

The main drawback of the Real Representation Theorem is that it is only
a \emph{sufficient} condition for membership because $f>0$ on $S(T)$
cannot be replaced by $f\ge 0$ (for example, a nonzero polynomial having a
zero in the interior of a polytope obviously never can allow Handelman's
representation adressed above). The criterion we will prove in Section
\ref{criterionsection}, Theorem
\ref{aboutissement} below, theoretically is necessary and sufficient. We say
``theoretically'' since it assumes the existence of a certain admissible
identity and there is a trivial identity (namely $f=1\cdot f$) that is
admissible if and only if $f\in T$. The Real Representation Theorem comes out
as a special case since another trivial identity (namely $f=f\cdot 1$) is
admissible if $f>0$ on $S(T)$. Our criterion yields new insights when
non-trivial admissible identities can be found. It is not of purely geometric
but also of arithmetic nature.

A slightly less general criterion for membership in \emph{preorderings}
has recently been proved by Scheiderer \cite[Proposition 3.10]{s3}. It has
been very successfully applied to partially extend Schm\"udgen's
representation from positive to nonnegative polynomials on compact basic
closed semialgebraic sets. Section \ref{preordersection} is devoted to the
question in how far our criterion goes beyond recent work of Scheiderer,
Kuhlmann, Marshall and Schwartz on preorderings. We will see that our
Theorem \ref{aboutissement} can easily be deduced from their work in the case
of preorderings, but in the general case, a central lemma in their proof
is no longer true (see Example \ref{exemplecontraire}). Our approach is
therefore not only different from theirs but also applies to a 
significantly more general situation.

In Section \ref{handelmansection}, we apply our criterion to give for the
first time a purely ring-theoretic proof of a nice
theorem of Handelman saying inter alia the following: If some odd power of a
real polynomial in several variables has only nonnegative coefficients, then
so do all sufficiently high powers. See Theorem \ref{handelman} and
Corollary \ref{odd}.

This example will nicely illustrate the following general principle:
Even if a semiring $T\subset A$ is not archimedean, there is always a
biggest subring $O_T(A)\subset A$ such that $T\cap O_T(A)$ is archimedean
This follows from the important Proposition \ref{archprop} below. So with
some additional difficulties
(namely determining $S(O_T(A))$), our membership criterion
also gives information about non-archimedean semirings.

\begin{proposition}\label{archprop}
Let $T$ be a semiring of $A$. Then
$$O_T(A):=\{a\in A\mid N\pm a\in T\text{\ for some $N\in\N$}\}$$
is a subring of $A$,
the \emph{ring of $T$-bounded elements} of $A$.
Moreover, $T$ is archimedean if and only if $O_T(A)=A$.
\end{proposition}

\begin{proof}
Obviously, $0,1\in O_T(A)$ since $0\pm 0=0\in T$ and
$1\pm 1\in\{0,2\}\subset T$. It is immediate from the definition
of $O_T(A)$ that $-O_T(A)\subset O_T(A)$. That $O_T(A)$ is closed under
addition, follows easily from $T+T\subset T$. To see that it is closed
under multiplication, use the two identities
$$3N^2\mp ab=(N\mp a)(N+b)+N(N\pm a)+N(N-b)$$
and that $T$ is closed under multiplication and addition. We leave the
second statement to the reader.
\end{proof}

Without going into details, we make some final remarks on the space $S(T)$.
There is a larger topological space
one could naturally associate to a semiring $T$ of a ring $A$, namely the
subspace $\Sper_T(A)$ of the so-called real spectrum $\Sper(A)$ of $A$
consisting of all so-called orderings of the ring $A$ lying over $T$
(see, e.g., \cite[4.1]{pd}).
Since $S(T)\subset\Sper_T(A)$ via a canonical embedding, all our results
will also be true for $\Sper_T(A)$. If $T$ is an archimedean semiring, then
$S(T)$ equals $(\Sper_T(A))^\text{max}$, the space of maximal orderings of
$A$ lying above $T$. When $T$ is not archimedean, $\Sper_T(A)$ is
certainly preferable to $S(T)$ (for example, $\Sper_T(A)$ is even then always
quasi-compact). However, we feel that in the context of archimedean
semirings we encounter here the usage of $\Sper_T(A)$ has only disadvantages.
For example, unlike $S(T)$, $\Sper_T(A)$ can usually not be really identified
with a concrete subset of $\R^n$. Confer also \cite[2.3]{s3}.

\section{The membership criterion}\label{criterionsection}

We begin by introducing some notation. For $\al\in\N^n$, we write
$$|\al|:=\al_1+\dots+\al_n,$$
so that the monomial
$$\x^\al:=X_1^{\al_1}\dotsm X_n^{\al_n}$$
has degree $|\al|$. For $x\in\R^n$, $\|x\|$ always denotes the $1$-norm
of $x$, i.e.,
$$\|x\|:=|x_1|+\dots+|x_n|.$$
Correspondingly,
$$B_r(x):=\{y\in\R^n\mid\|y-x\|<r\}\qquad(x\in\R^n,0<r\in\R)$$
denotes the open ball around $x$ of radius $r$ with respect to the $1$-norm and
$$\overline{B_r(x)}=\{y\in\R^n\mid\|y-x\|\le r\}$$
its closure. Like all norms, the $1$-norm defines the usual topology
on $\R^n$. The reason for our choice of this norm is that
$\|\al\|=|\al|$ for $\al\in\N^n$. Despite this equality, we want to keep
both notations since $|\al|=k$ will mean implicitly $\al\in\N^n$
(and that $\al$ plays the role of a tuple of exponents of a monomial
$\x^\al$).
We introduce the compact set
\begin{align*}
\De&:=(\overline{B_1(0)}\setminus B_1(0))\cap\nn^n=V(\{X_1+\dots+X_n-1\})
\cap\nn^n\\
&=\{x\in\nn^n\mid \|x\|=1\}\subset\R^n.
\end{align*}
For a given set $P\subset\rx$, we denote by $P^+$
its subset of all polynomials which have only nonnegative coefficients
and by $P^\ast$ its subset of all homogeneous polynomials (i.e., all of
whose nonzero monomials have the same degree).

Descartes had already the idea to relate the geometric properties of a real
polynomial directly to combinatorial properties of the family of signs of its
coefficients. His law of signs says that a real polynomial in one variable has
not more positive real roots than it has sign changes in the sequence of its
coefficients, and the difference is even \cite[Theorem 2.34]{bpr}. Given a
sequence of signs, a good guess for the number of positive real roots of a
corresponding polynomial would therefore perhaps be the number of these
sign changes. Viro extended this naive rule of guessing the topological shape
of the real zero set of a polynomial to the case of several variables.
Given a pattern of signs, he can construct a corresponding polynomial
whose real zero set has exactly the guessed shape. This is Viro's method for
constructing real hypersurfaces with prescribed topology \cite{vir}.

The starting point for the proof of our criterion is yet another idea in this
vein going back to P\'olya. Suppose $f\in\rx^\ast$. P\'olya relates the
geometric behaviour of $f$ on the nonnegative orthant $[0,\infty)^n$ with
the signs of the coefficients of a ``refinement'' of $f$. Due to homogeneity,
$f$ can just as well be looked at on $\De$ instead of $[0,\infty)^n$.
Multiplying $f$ by $X_1+\dots+X_n$ does not change $f$ on $\De$ but
``refines'' the pattern of signs of its coefficients. When we repeat this
multiplication sufficiently often, it turns out that the obtained pattern
reflects more and more the geometric sign behaviour of $f$ on $[0,\infty)^n$.
The exact statement we will need is formulated in Lemma \ref{localpolya}
below. Whereas previous works of the author \cite{sw1}\cite{sw3}\cite{sw4}
(see Remark \ref{quark} below) required only P\'olya's original theorem,
we need this time really a more local version where we look at $f$ only on
a closed subset $U$ of $\De$. Nevertheless, the proof goes exactly along the
lines of P\'olya (cf. \cite{pol}\cite{pr}). We include it for the
convenience of the reader.

\begin{lemma}\label{localpolya}
Suppose $f\in\rx^\ast$ has degree $d$ and $U\subset\De$ is
closed such that $f>0$ on $U$. Then there is $k_0\in\N$ such that for
all $k\ge k_0$ and $\al\in\N^n$ with ($k+d\neq 0$ and) $\frac\al{k+d}\in U$,
the coefficient of $\x^\al$ in $(X_1+\dots+X_n)^kf$ is positive.
\end{lemma}

\begin{proof}
Write $f=\sum_{|\be|=d}a_\be\x^\be$, $a_\be\in\R$. We know that
$$(X_1+\dots+X_n)^k=\sum_{|\ga|=k}\frac{k!}{\ga_1!\dotsm\ga_n!}\x^\ga$$
for $k\in\N$. Of course, if $\al\in\N^n$ with $\frac\al{k+d}\in U\subset\De$,
then $|\al|=k+d$.
Now for any $\al\in\N^n$ with $|\al|=k+d$, the coefficient of
$\x^\al$ in $(X_1+\dots+X_n)^kf$ equals
\begin{align*}
&\sum_{\genfrac{}{}{0pt}{}{|\be|=d,|\ga|=k}{\be+\ga=\al}}
\frac{k!}{\ga_1!\dotsm\ga_n!}a_\be\\
=&\sum_{\genfrac{}{}{0pt}{}{|\be|=d,|\ga|=k}{\be+\ga=\al}}
 \frac{k!}{(\al_1-\be_1)!\dotsm(\al_n-\be_n)!}a_\be
 \qquad\text{(terms of the sum do not depend on $\ga$)}\\
=&\sum_{\genfrac{}{}{0pt}{}{|\be|=d}{\be\le\al}}
 \frac{k!}{(\al_1-\be_1)!\dotsm(\al_n-\be_n)!}a_\be
 \qquad\text{($\be\le\al$ understood componentwise)}\\
=&\frac{k!(k+d)^d}{\al_1!\dotsm\al_n!}
 \sum_{\genfrac{}{}{0pt}{}{|\be|=d}{\be\le\al}}
 a_\be\prod_{i=1}^n\frac{\al_i!}{(\al_i-\be_i)!(k+d)^{\be_i}}
 \qquad\text{(using $|\be|=d$, provided $k+d\neq 0$)}\\
=&\frac{k!(k+d)^d}{\al_1!\dotsm\al_n!}\sum_{|\be|=d}
 a_\be\prod_{i=1}^n
 \left(\frac{\al_i}{k+d}\right)^{\be_i}_{\frac 1{k+d}}
 \qquad\text{(abbreviating $(a)_b^m:=\prod_{i=0}^{m-1}(a-ib)$).}
\end{align*}
Note that $(a)_0^m=a^m$ to understand the idea behind the notation $(a)_b^m$
just introduced. Also note that the condition $\be\le\al$ has been dropped in
the index of summation in the last expression. This is justified since all the
corresponding additional terms in the sum are zero. Now we see that the
coefficient in question equals (assuming $k+d\neq 0$) up to a positive factor
$$f_{\frac 1{k+d}}\left(\frac\al{k+d}\right)$$
where we define
$$f_\ep:=\sum_{|\be|=d}a_\be(X_1)^{\be_1}_\ep\dotsm(X_n)^{\be_n}_\ep\in\rx$$
for all $\ep\in\nn$. Obviously, $f_\ep$ converges to $f_0=f$ uniformly on $U$
when $\ep\to 0$. Since $U$ is compact and $f>0$ on $U$, there is
$k_0\in\N$ such that $f_{\frac 1{k+d}}>0$ on $U$ for all $k\ge k_0$, in
particular $$f_{\frac 1{k+d}}\left(\frac\al{k+d}\right)>0$$ whenever
($k+d\neq 0$ and) $\frac\al{k+d}\in U$.
\end{proof}

We draw from this P\'olya's theorem as a corollary although we will never use
it later. Note that P\'olya's theorem follows as easily by taking independently
of $x\in\De$ the same identity $f=f\cdot 1$ in condition (\ref{fhid}) of
Lemma \ref{nonnelya} below.

\begin{corollary}[P\'olya]\label{polya}
Suppose $f\in\rx^\ast$ and $f>0$ on $\De$. Then
$$(X_1+\dots+X_n)^kf\in\rx^+$$ for large $k\in\N$.
\end{corollary}

\begin{proof}
Set $U=\De$ in Lemma \ref{localpolya}.
\end{proof}

It is perhaps worth pointing out that P\'olya's theorem is closely related to
Bernstein polynomials. See \cite[Theorem 1.3]{far} for a theorem on
(generalized) Bernstein polynomials (in several variables) which is nothing
else than a version of P\'olya's theorem. Via this connection, P\'olya's
Theorem for the case of two variables (i.e., when $n=2$) is connected to
Descartes' law of signs mentioned above \cite[Section 10.27]{bpr}.
For technical reasons, it is very convenient to have the following 
evident consequence of Lemma \ref{localpolya} available.

\begin{lemma}\label{papote}
Suppose $f\in\rx^\ast$ and $U\subset\De$ is
closed such that $f>0$ on $U$. Then there is $k_0\in\N$ such that for
all $k\ge k_0$ and $0\neq\al\in\N^n$ with $\frac\al{|\al|}\in U$, the
coefficient of $\x^\al$ in $(X_1+\dots+X_n)^kf$ is nonnegative.
\end{lemma}

\begin{proof} Without loss of generality $f\neq 0$. Set $d:=\deg f$
and choose $k_0$ like in the previous lemma. Let $k\ge k_0$
and $0\neq\al\in\N^n$ with $\frac\al{|\al|}\in U$. If $|\al|=k+d$, $\x^\al$
has a positive coefficient in $(X_1+\dots+X_n)^kf$ by the choice of $k_0$.
If $|\al|\neq k+d$, the coefficient
of $\x^\al$ in this same polynomial is zero since it is a homogeneous
polynomial of degree $k+d$.
\end{proof}

The next lemma reminds already a bit of Theorem \ref{aboutissement} below.
But note that the $g_i$ and $h_i$ are allowed to depend on $x$. The idea is
to apply P\'olya's refinement process locally on the $g_i$ while the $h_i$ do
not disturb too much. Note that we do no longer assume that $f$ is homogeneous.
Also observe that the hypotheses imply $f\ge 0$ on $\De$.

\begin{lemma}\label{nonnelya}
Let $f\in\rx$. Suppose that for every $x\in\De$ there are $m\in\N$,
$g_1,\dots,g_m\in\rx^\ast$ and $h_1,\dots,h_m\in\rx^+$ such that
\begin{enumerate}[(a)]
\item $f=g_1h_1+\dots+g_mh_m$ and\label{fhid}
\item $g_1(x)>0,\dots,g_m(x)>0$.\label{gpos}
\end{enumerate}
Then there exists $k\in\N$ such that $(X_1+\dots+X_n)^kf\in\rx^+$.
\end{lemma}

\begin{proof}
Choose a family $(\ep_x)_{x\in\De}$ of real numbers $\ep_x>0$ such that
for every $x\in\De$, there are $m\in\N$, $g_1,\dots,g_m\in\rx^\ast$ and
$h_1,\dots,h_m\in\rx^+$
satisfying (\ref{fhid}) and not only (\ref{gpos}) but even
\begin{equation}\label{gposplus}
g_i>0\text{\ on\ }\overline{B_{2\ep_x}(x)}\cap\De\qquad
\text{for $i\in\{1,\dots,m\}$.}
\end{equation}
The family $(B_{\ep_x}(x))_{x\in\De}$ is an
open covering of $\De$. Since $\De$ is compact, there is a finite
subcovering, i.e., a finite set $D\subset\De$ for which
$\De\subset\bigcup_{x\in D}B_{\ep_x}(x)$, in particular
$$\De=\bigcup_{x\in D}(\overline{B_{\ep_x}(x)}\cap\De).$$
As $D$ is finite, it suffices to show for fixed $x\in D$, that there is
$k_0\in\N$ such that for all $k\ge k_0$ and all $0\neq\al\in\N^n$ with
\begin{equation}\label{toshow}
\frac\al{|\al|}\in\overline{B_{\ep_x}(x)},
\end{equation} the coefficient of $\x^\al$ in $(X_1+\dots+X_n)^kf$ is
nonnegative (note that $\frac\al{|\al|}\in\De$ is automatic).

Therefore fix $x\in D$. By choice of $\ep_x$, we find $m\in\N$,
$g_1,\dots,g_m\in\rx^\ast$ and $h_1,\dots,h_m\in\rx^+$ satisfying
(\ref{fhid}) and (\ref{gposplus}). For every $i\in\{1,\dots,m\}$, the
positivity condition (\ref{gposplus}) enables us to apply Lemma
\ref{papote} to $g_i$, yielding $k_i\in\N$ such that for all
$k\ge k_i$ and all $0\neq\be\in\N^n$ with
\begin{equation}\label{betacond}
\frac\be{|\be|}\in\overline{B_{2\ep_x}(x)},
\end{equation}
the coefficient of $\x^\be$ in $(X_1+\dots+X_n)^kg_i$ is nonnegative
(use that $\frac\be{|\be|}\in\De$ is automatic).
Choose moreover $1\le l\in\N$ so large that
\begin{equation}\label{lchoice}
\frac{2|\ga|}l\le\ep_x
\end{equation}
for all $\ga\in\N^n$ for which the coefficient of $\x^\ga$ in at least
one of the polynomials $h_1,\dots,h_m$ does not vanish. Set
$$k_0:=\max\{k_1,\dots,k_n,l\}.$$
Let $k\ge k_0$ and suppose $0\neq\al\in\N^n$ satisfies $(\ref{toshow})$. Fix
$i\in\{1,\dots,m\}$. By equation (\ref{fhid}), it is enough to show that the
coefficient of $\x^\al$ in $(X_1+\dots+X_n)^kg_ih_i$ is nonnegative. This
coefficient is of course a sum of certain products of coefficients of
$(X_1+\dots+X_n)^kg_i$ and $h_i$. But all the concerned products are
nonnegative. Indeed, consider $\be,\ga\in\N^n$ with $\be+\ga=\al$
(i.e., $\x^\be\x^\ga=\x^\al$) such that the corresponding coefficients of
$\x^\be$ in $(X_1+\dots+X_n)^kg_i$ and $\x^\ga$ in $h_i$ do not vanish.
The latter coefficient is positive since $h_i\in\rx^+$. We show that the other
one is positive, too. From degree consideration it is trivial that
$|\be|\ge k\ge k_0\ge l\ge 1$ which implies together with the now
satisfied condition (\ref{lchoice})
\begin{equation}\label{begarel}
\frac{2|\ga|}{|\be|}\le\ep_x.
\end{equation}
We exploit this to verify condition (\ref{betacond}) which is all we
need since $k\ge k_0\ge k_i$:
\begin{align*}
\left\|\frac\be{|\be|}-x\right\|&\le
\left\|\frac{\be}{|\be|}-\frac\al{|\al|}\right\|+
\underbrace{\left\|\frac\al{|\al|}-x\right\|}_
{\le\ep_x\text{\ by (\ref{toshow})}}
\le\ep_x+\left\|\frac{|\al|\be-|\be|\al}
  {|\al||\be|}\right\|\\
&=\ep_x+\frac 1{|\al||\be|}
\|\underbrace{|\al|(\ga-\al)-\overbrace{|\ga-\al|}^{=|\ga|-|\al|}\al}_
{=|\al|\ga-|\al|\al-|\ga|\al+|\al|\al}\|\\
&=\ep_x+\frac{\||\al|\ga-|\ga|\al\|}{|\al||\be|}
\le\ep_x+\frac{\||\al|\ga\|+\||\ga|\al\|}{|\al||\be|}\\
&=\ep_x+\frac{2|\al||\ga|}{|\al||\be|}=\ep_x+\frac{2|\ga|}{|\be|}
\overset{(\ref{begarel})}\le 2\ep_x
\end{align*}
\end{proof}

Now we deal with the case where the $g_i$ are no longer assumed to be
homogeneous.

\begin{lemma}\label{modsingle}
Let $f\in\zx$ such that for all $x\in\De$,
there exist $m\in\N$, $g_1,\dots,g_m\in\zx$ and $h_1,\dots,h_m\in\zx^+$
such that
\begin{enumerate}[(a)]
\item $f=g_1h_1+\dots+g_mh_m$ and\label{fgh}
\item $g_1(x)>0,\dots,g_m(x)>0$.\label{ggg}
\end{enumerate}
Then $f$ is modulo the principal ideal $\Z[\x](X_1+\dots+X_n-1)$
congruent to a polynomial without negative coefficients.
\end{lemma}

\begin{proof}
For every $x\in \De$, choose
$m_x\in\N$, $g_{x1},\dots,g_{x{m_x}}\in\zx$ and
$0\neq h_{x1},\dots,h_{x{m_x}}\in\zx^+$ according
to (\ref{fgh}) and (\ref{ggg}). Setting
\begin{equation}\label{uxdef}
U_x:=\{y\in\De\mid g_{x1}(y)>0,\dots,g_{x{m_x}}(y)>0\},
\end{equation}
we have $x\in U_x$ for $x\in\De$. Therefore
$(U_x)_{x\in\De}$ is an open covering of the compact set $\De$ and
possesses a finite subcovering, i.e., there is a finite set $D\subset\De$
such that
\begin{equation}\label{cover}
\De=\bigcup_{x\in D}U_x.
\end{equation}
Choose an upper bound $d\in\N$ for the
degrees of the (in each case $m_x$) terms appearing in the sums on the right
hand sides of
the equations (\ref{fgh}) corresponding to the finitely many $x\in D$,
i.e.,
$$d\ge\deg g_{xi}+\deg h_{xi}\qquad\text{for all $x\in D$ and
  $i\in\{1,\dots,m_x\}$}.$$
Fix for the moment such a pair $(x,i)$ and choose $d',d''\in\N$ such that
$d=d'+d''$, $d'\ge\deg g_{xi}$ and $d''\ge\deg h_{xi}$. Write
$g_{xi}=\sum_{k=0}^{d'}p_k$ and $h_{xi}=\sum_{k=0}^{d''}q_k$ where
$p_k,q_k\in\zx$ are homogeneous of degree $k$ (if not zero). Set
$$g_{xi}':=\sum_{k=0}^{d'}(X_1+\dots+X_n)^{d'-k}p_k\quad\text{and}\qquad
  h_{xi}':=\sum_{k=0}^{d''}(X_1+\dots+X_n)^{d''-k}q_k.$$
Now $g_{xi}'$ and $h_{xi}'$ are homogeneous polynomials whose product
is (homogeneous) of degree $d$ (if not zero). Then
$g_{xi}'\equiv g_{xi}$ and $h_{xi}'\equiv h_{xi}$ modulo
$\Z[\x](X_1+\dots+X_n-1)$, in particular, $g_{xi}'$ coincides with $g_{xi}$
on $\De$. Moreover, $h_{xi}'$ inherits the property of having no negative
coefficients from $h_{xi}$. For every $x\in D$,
\begin{equation}\label{fprimedef}
f_x':=g_{x1}'h_{x1}'+\dots+g_{xm_x}'h_{xm_x}'\in\zx^\ast
\end{equation}
is homogeneous of degree $d$ (unless zero) and congruent to $f$ modulo
$\Z[\x](X_1+\dots+X_n-1)$. For $x,y\in D$, $f_x'-f_y'$ is therefore homogeneous
and at the same time a multiple of $X_1+\dots+X_n-1$. Hence actually
$f_x'=f_y'$, i.e., there is $f'\in\zx$ such that $f'=f_x'$ for all $x\in D$
and $f'\equiv f$ modulo $\Z[\x](X_1+\dots+X_n-1)$.

We want to apply Lemma
\ref{nonnelya} to $f'$. The hypotheses are now rather easy to verify:
Let $x\in\De$. By (\ref{cover}), we find $x\in D$ such that $x\in U_x$.
Set $m:=m_x$, $g_i:=g_{xi}'$ and $h_i:=h_{xi}'$ for $i\in\{1,\dots,m\}$.
Then equation (\ref{fprimedef}) becomes condition $(\ref{fhid})$ in
Lemma \ref{nonnelya} (with $f'$ instead of $f$). To verify (\ref{gpos}) of
Lemma \ref{nonnelya}, use that $g_i=g_{xi}'$ equals $g_{xi}$ on $\De$ which is
positive in $x\in U_x\subset\De$ by (\ref{uxdef}).
By Lemma \ref{nonnelya}, we get therefore $k\in\N$ such that
$(X_1+\dots+X_n)^kf'$ has no negative coefficients. But this polynomial
is congruent to $f'$ which is in turn congruent to $f$ modulo
$\Z[\x](X_1+\dots+X_n-1)$.
\end{proof}

Compared to the lemma we just proved, the next statement has the big
advantage that the principal ideal $\Z[\x](X_1+\dots+X_n-1)$ can be replaced by
any larger ideal. On the other hand, the $h_i$ are no longer allowed to
depend on $x$. This disadvantage is made more tolerable by the fact that
only those $x$ have to be considered where $f$ vanishes. In the previous
(but not in the next) lemma this fact is implicitly obvious since
$f=f\cdot 1$ is an admissible identity at the points where $f$ is positive.

\begin{lemma}\label{modmulti}
Let $I$ be an ideal of $\zx$ such that $X_1+\dots+X_n-1\in I$. Suppose
$m\in\N$, $f\in\zx$ and $h_1,\dots,h_m\in\zx^+$ such that
\begin{enumerate}[(a)]
\item $f\ge 0$ on $V(I)\cap\nn^n$ and\label{nnc}
\item for all $x\in V(I\cup\{f\})\cap\nn^n$, there exist\label{ccc}
$g_1,\dots,g_m\in\zx$ such that
\begin{enumerate}
\item[(i)] $f=g_1h_1+\dots+g_mh_m$ and
\item[(ii)] $g_1(x)>0,\dots,g_m(x)>0$.
\end{enumerate}
\end{enumerate}
Then $f$ is modulo $I$ congruent to a polynomial without negative coefficients.
\end{lemma}

\begin{proof}
Set $U:=\{x\in\De\mid f(x)>0\}$ and introduce the set
$W\subset\De$ of all $x\in\De$ for which there are $g_1,\dots,g_m$
fulfilling (i) and (ii). The sets $U$ and $W$ are open in $\De$ and
\begin{equation}\label{towards}
V(I)\cap\nn^n\subset U\cup W
\end{equation}
by (\ref{nnc}) and (\ref{ccc}).
By Hilbert's Basis Theorem, every ideal of $\zx$ is finitely generated. In
particular, we find $s\in\N$ and $p_1,\dots,p_s\in\zx$ such that
$$I=\Z[\x]p_1+\dots+\Z[\x]p_s+\Z[\x](X_1+\dots+X_n-1).$$
Setting $p:=\sum_{i=1}^sp_i^2\in I$, we have $p\in I$, $p\ge 0$ on $\R^n$ and
\begin{equation}\label{stayaway}
\text{$p>\ep$ on $\De\setminus (U\cup W)$ \qquad for some $\ep>0$.}
\end{equation}
The latter follows from $p>0$ on $\De\setminus V(I)$, (\ref{towards})
and the compactness of $\De\setminus (U\cup W)$.

Now we distinguish two cases. First case: $W=\emptyset$. From (\ref{stayaway})
and the boundedness of $f$ on the compact set $\De\setminus U$, we get
$k\in\N$ such that $f':=f+kp>0$ on $\De\setminus U$.
On the other hand, $f'=f+kp\ge f>0$ on $U$. Altogether we get $f'>0$ on $\De$.
Now we can clearly apply Lemma \ref{modsingle} to $f'$. In fact, for every
$x\in\De$, $f'=f'\cdot 1$ serves as an identity as required in (\ref{fgh}) of
that lemma. Hence that lemma yields that $f'$ is
congruent to a polynomial without negative coefficients modulo
$\Z[\x](X_1+\dots+X_n-1)\subset I$. But $f\equiv f+kp=f'$
modulo $I$.

Second case: $W\neq\emptyset$. All we really use from $W\neq\emptyset$ is
that $f\in\Z[\x]h_1+\dots+\Z[\x]h_m$ by (i), i.e., we find
$q_1,\dots,q_m\in\zx$ such that
\begin{equation}\label{idealmember}
f=q_1h_1+\dots+q_mh_m.
\end{equation}
From (\ref{stayaway}) and the boundedness of $q_1,\dots,q_m$ on the compact set
$\De\setminus (U\cup W)$, it follows that we can choose $k\in\N$ such that
\begin{equation}\label{fpd}
g_i^{(0)}:=q_i+kp>0\qquad\text{on $\De\setminus (U\cup W)$ for all
$i\in\{1,\dots,m\}$.}
\end{equation}
We will apply Lemma \ref{modsingle} to
\begin{equation}\label{defaultidentity}
f':=g_1^{(0)}h_1+\dots+g_m^{(0)}h_m.
\end{equation}
Note that
\begin{equation}\label{expl}
f'\overset{(\ref{fpd})}=
\underbrace{q_1h_1+\dots+q_mh_m}_{\text{$=f$ by (\ref{idealmember})}}+
kp(\underbrace{h_1+\dots+h_m}_{\ge 0\text{\ on $\nn^n$}})\ge f
\qquad\text{on $\nn^n$.}
\end{equation}
To check its applicability, let
$x\in\De$. We consider three different subcases:

First, consider the case where $x\in U$. Then $f'(x)\ge f(x)>0$ and
\begin{equation}\label{trividentity}
f'=f'\cdot 1
\end{equation}
is an identity as demanded in (\ref{fgh}) of Lemma
\ref{modsingle}.

Second, suppose $x\in W$. By definition of $W$, we can choose
$g_1,\dots,g_m\in\zx$ satisfying $(i)$ and $(ii)$.
Set $g_i':=g_i+kp$ for $i\in\{1,\dots,m\}$. Then
\begin{equation}\label{nondefaultidentity}
\begin{split}
f'&\overset{(\ref{expl})}=f+kp(h_1+\dots+h_m)\\
&\overset{(i)}=g_1h_1+\dots+g_mh_m+kp(h_1+\dots+h_m)=g_1'h_1+\dots+g_m'h_m
\end{split}
\end{equation}
serves as a relation as required in (\ref{fgh}) of Lemma \ref{modsingle}. Note
that
$$g_i'(x)=g_i(x)+kp(x)\ge g_i(x)\overset{(ii)}>0\qquad
  \text{for $i\in\{1,\dots,m\}$.}$$

Third and last, for all $x\in\De\setminus(U\cup W)$, (\ref{fpd})
allows us to use one and the same equation for (\ref{fgh}) of Lemma
\ref{modsingle}, namely (\ref{defaultidentity}).

All in all, Lemma \ref{modsingle} applies now to $f'$, i.e., $f'$ is
congruent to a polynomial without nonnegative coefficients modulo
$\Z[\x](X_1+\dots+X_n-1)\subset I$. But $f\equiv f+kp=f'$ modulo $I$.
\end{proof}

\begin{remark}\label{perte}
In Lemma \ref{modsingle}, the $h_i$ are permitted to depend on
$x$. In the proof of Lemma \ref{modmulti}, we do not exploit this too much.
Indeed, the three used identities (\ref{trividentity}),
(\ref{nondefaultidentity}) and (\ref{defaultidentity}) are based on the
same $h_i$ except (\ref{trividentity}) which is a trivial identity.
\end{remark}

For any element $a\in A$, we set
$$S_{a=0}(T):=\{x\in S(T)\mid a(x)=0\}.$$
Now we attack the main theorem. Note that its hypotheses imply that all
$t_i$ vanish
on $S_{a=0}(T)$.

\begin{theorem}\label{aboutissement}
Let $T$ be a weakly divisible archimedean semiring of $A$ and $a\in A$. Suppose
$a\ge 0$ on $S(T)$ and there is an identity $a=b_1t_1+\dots+b_mt_m$ with
$b_i\in A$, $t_i\in T$ such that $b_i>0$ on $S_{a=0}(T)$ for all
$i\in\{1,\dots,m\}$. Then $a\in T$.
\end{theorem}

\begin{proof}
If the ring homomorphism $\Z\to A$ is not injective, then $-1\in T$ whence
$T=\Z+T=A$. Therefore we assume from now on that $A$ contains $\Z[\frac 1r]$
as a subring and $\frac 1r\in T$ for some integer $r\ge 2$. Because $T$ is
archimedean, we find for every $c\in A$ some $N_c\in\N$ with
$N_c\pm c\in T$. The topological space
$$S:=\prod_{c\in A}[-N_c,N_c]$$
is compact by Tychonoff's theorem. From the hypotheses of the theorem, it
follows that a certain intersection of closed subsets of $S$ is empty:
\begin{align}
\bigcap_{c,d\in A}&\{\ph\in S\mid\ph(c)+\ph(d)-\ph(c+d)=0\}&\cap\label{iplus}\\
\bigcap_{c,d\in A}&\{\ph\in S\mid\ph(c)\ph(d)-\ph(cd)=0\}&\cap\label{itimes}\\
&\{\ph\in S\mid\ph(1)=1\}&\cap\label{ione}\\
\bigcap_{t\in T}&\{\ph\in S\mid\ph(t)\ge 0\}&\cap\label{it}\\
&\{\ph\in S\mid\ph(a)\le 0\}&\cap\label{if}\\
\bigcup_{i=1}^m&\{\ph\in S\mid\ph(b_i)\le 0\}&&=\emptyset\label{ia}
\end{align}
All sets appearing as subexpressions of (\ref{iplus})--(\ref{ia}) are
closed. This is easy to see: Use that $\{0\},\{1\},\nn,(-\infty,0]$ are
closed subsets of $\R$, that the projection maps $S\to\R:\ph\mapsto\ph(c)$
($c\in A$) are continuous (the characteristic property of the product
topology), that $+,-,\cdot:\R\times\R\to\R$ are continuous and that finite
unions and arbitrary intersections of closed sets are again closed.

Since $S$ is compact, some finite subintersection of (\ref{iplus})--(\ref{ia})
is already empty. In particular, (\ref{iplus})--(\ref{ia}) is already empty
if the intersection in (\ref{iplus}) and in (\ref{itimes}) runs only over
certain finitely many $c,d\in A$. Let $\bar y=(y_1,\dots,y_n)$ be the
collection of $\frac 1r$, all $b_i$, $t_i$ and these $c,d$. We claim that all
hypotheses of the theorem remain valid for $(\Z[\bar y],T\cap\Z[\bar y])$
instead of $(A,T)$.

Indeed, first of all, $T\cap\Z[\bar y]$ inherits the property
of being a weakly divisible archimedean semiring from $T$.
Second, the identity from the hypotheses remains trivially satisfied
(do not forget that $\Z[\bar y]$ contains $a$ since it contains
all $b_i,t_i$ and it is a ring). Third and last, it remains to check that the
geometric hypotheses stay valid.
To this purpose, let $\ph:\Z[\bar y]\to\R$ be a ring
homomorphism with $\ph(T\cap\Z[\bar y])\subset\nn$. Extend $\ph$ to a
map
$$\ps:A\to\R:c\mapsto\begin{cases}
                      0&\text{if $c\notin\Z[\bar y]$,}\\
                      \ph(c)&\text{if $c\in\Z[\bar y].$}
                      \end{cases}
$$
We have $\ps\in S$: If $c\in A\setminus\Z[\bar y]$, then
$\ps(c)=0\in[-N_c,N_c]$. If $c\in\Z[\bar y]$, then
$N_c\pm c\in T\cap\Z[\bar y]$ by choice of $N_c$, whence
\begin{equation}\label{malheure}
N_c\pm\ps(c)=N_c\pm\ph(c)=\ph(N_c\pm c)\in\ph(T\cap\Z[\bar y])
\subset\nn
\end{equation}
showing also in this case $\ps(c)\in[-N_c,N_c]$. Since $\ph$ is a ring
homomorphism, its extension $\ps$ satisfies the corresponding homomorphy
conditions on $\bar y$. Therefore $\ps\in S$ lies in the finite
subintersection of (\ref{iplus})--(\ref{itimes}) that led us above to
the choice of $\bar y$. It is even easier to see that $\ps$ also
lies in the intersection (\ref{ione})--(\ref{it}). Because $\ps$ cannot
lie in the empty set, $\ps$ cannot lie in both (\ref{if}) and (\ref{ia}).
Hence $\ph(a)=\ps(a)>0$ or, for all $i$, $\ph(b_i)=\ps(b_i)>0$.
In the latter case $\ph(a)=\ph(b_1)\ph(t_1)+\dots+\ph(b_m)\ph(t_m)\ge 0$.
Altogether, this shows $a\ge 0$ on $S(T\cap\Z[\bar y])$ and
$b_i>0$ on $S_{a=0}(T\cap\Z[\bar y])$ for all $i$.
In other words, the hypotheses of the theorem are valid for
$\left(\Z\left[\bar y\right],T\cap\Z\left[\bar y\right]\right)$
instead of $(A,T)$. So we can assume from now on that
\begin{equation}\label{fingen}
A=\Z[\bar y],
\end{equation}
i.e., that $A$ is finitely generated as a ring. We still let $r\ge 2$ be an
integer such that $\frac 1r\in T$. We assume that
\begin{equation}\label{xint}
y_i\in T\qquad\text{for all $i\in\{1,\dots,n\}.$}
\end{equation}
This is justified by the fact that $T$ is archimedean since there is $N\in\N$
with $N+y_i\in T$ and we may replace $y_i$ by $N+y_i$
(this does not affect (\ref{fingen})). Moreover, the
assumption
\begin{equation}\label{sumsup}
y_1+\dots+y_n=1
\end{equation}
is without loss of generality:
We can extend $\bar y$ by $N-(y_1+\dots+y_n)$ for any $N\in\N$ without
harming (\ref{fingen}). If we choose $N\in\N$ so large that
$N-(y_1+\dots+y_n)\in T$, then (\ref{xint}) remains valid at the same time.
Choosing this $N$ even more carefully, namely as a power $r$, establishes
(\ref{sumsup}) with a power of $r$ instead of $1$ on the right hand side.
Finally, divide each $y_i$ by this power of $r$ (cf. (\ref{fingen})).

Now consider the ring epimorphism $\zx\to\Z[\bar y]$ mapping
$X_i$ to $y_i$ for every $i\in\{1,\dots,n\}$. Calling its kernel $I$,
it induces a ring isomorphism $\zx/I\to A$ mapping $X_i+I$ to $y_i$.
Without loss of generality, we may assume
\begin{equation}\label{polrep}
A=\zx/I\qquad\text{and}\qquad y_i=X_i+I\text{\ for $i\in\{1,\dots,n\}$.}
\end{equation}
As explained in Section \ref{semiringsection}, we then have the concrete
description (\ref{concret}) of $S(T)$, i.e.,
\begin{equation}\label{repetition}
S(T)=\{x\in V(I)\mid\text{$t(x)\ge 0$ for all $t\in T$}\}\subset\R^n.
\end{equation}
The geometric part of the hypotheses of our theorem implies
\begin{align}
\bigcap_{t\in T}&\{x\in \De\cap V(I)\mid t(x)\ge 0\}&\cap\label{is1}\\
&\{x\in \De\cap V(I)\mid a(x)\le 0\}&\cap\label{is2}\\
\bigcup_{i=1}^m&\{x\in \De\cap V(I)\mid b_i(x)\le 0\}&&=\emptyset.\label{is3}
\end{align}
This is analogous to the above intersection (\ref{iplus})--(\ref{ia}): $V(I)$
plays the role of subintersection (\ref{iplus})--(\ref{ione}) and $\De$ plays
the role of $S$. Exactly as we would even get an empty intersection 
 in (\ref{iplus})--(\ref{ia}) above with $S$ replaced by $\R^A$, we could
replace here $\De$ by $\R^n$. But in order to have an intersection of closed
subsets of a \emph{compact} space, we have defined all sets as subsets of
$S$ above and define them as subsets of $\De$ here. The fact that everything
now happens in $\R^n$ instead of $\R^A$ is important. It will allow us
to pass over to a finitely generated semiring $T'\subset T$. See Remark
\ref{passage} below.

As already pointed out, (\ref{is1})--(\ref{is3}) is an empty intersection of
closed sets in the compact space $\De\cap V(I)$. Hence it has a finite empty
subintersection. In particular, it is already empty if the intersection in
(\ref{is1}) runs only over finitely many (instead of all) $t\in T$. Let
$T'\subset A$ be the semiring generated by these finitely many $t$ and
$\frac 1r,y_1,\dots,y_n,t_1,\dots,t_m$. 

The semiring $T'$ still is (weakly divisible and) archimedean. According to
Proposition \ref{archprop} and (\ref{fingen}), this can be verified by checking
$y_1,\dots,y_n\in O_{T'}(A)$.
But this is immediate from $1+y_i\in 1+T'\subset T'$ and
$$1-y_i\overset{(\ref{sumsup})}=\sum_{j\neq i}y_j\in T'$$
for $i\in\{1,\dots,n\}$.

Next, we claim that $a\ge 0$ on $S(T')$ and $b_i>0$ on $S_{a=0}(T')$.
So let $x\in S(T')$. With respect to a description of $S(T')$ analogous
to (\ref{repetition}), we have of course $x\in V(I)$.
From (\ref{xint}) with $T'$ instead of $T$ and
(\ref{sumsup}), we obtain $x\in\De$.
Therefore $x$ is contained in the finite subintersection
of (\ref{is1}) which led to the choice of $T'$. So it cannot be contained
in intersection (\ref{is2})--(\ref{is3}). So, if $a(x)\le 0$, then
$b_i(x)>0$ for all $i$, whence $a(x)\ge 0$ (so actually $a(x)=0$)
by the identity from the hypotheses (recall that $t_1,\dots,t_m\in T'$).

Now, we see that the hypotheses of the theorem remain satisfied with
$T$ substituted by $T'$. As $T'\subset T$, it is consequently enough to
show the theorem for $T'$ instead of $T$. The advantage is that $T'$
is finitely generated as a semiring. For ease of notation, we work
again with $T$ instead of $T'$ but can assume from now on that
$T$ is a finitely generated semiring. But then we see that we could have
chosen $y_1,\dots,y_n$ fulfilling (\ref{fingen}) in such a way that they
generate $T$. Let us assume henceforth that we did so. Then it follows from
(\ref{polrep}) that
\begin{equation}\label{hasno}
T=\{p+I\mid p\in\zx^+\}.
\end{equation}
We see from this that
\begin{equation}\label{nnnnn}
S(T)=V(I)\cap\nn^n.
\end{equation}
Choose $g_1,\dots,g_m\in\zx$ and $h_1,\dots,h_m\in\zx^+$ such that $b_i=g_i+I$,
$t_i=h_i+I$ for all $i$ (use (\ref{hasno})). Now set
\begin{equation}\label{lacleman}
f:=g_1h_1+\dots+g_mh_m\in\zx
\end{equation}
which is nothing else than condition (i) in Lemma \ref{modmulti}.
The remaining hypotheses of Lemma \ref{modmulti} are now
provided by (\ref{sumsup}), (\ref{polrep}) and (\ref{nnnnn}). That
lemma yields that $f$ is congruent to a polynomial without negative
coefficients modulo $I$. By (\ref{hasno}), this means that $a=f+I\in T$.
\end{proof}

Together with Remark \ref{perte}, the next remark will tell us that the
intermediate results in this section have not been exploited to their full
extent. This gives hope that the just proved theorem can still be improved
at least in certain special situations.

\begin{remark}
In condition (\ref{ccc}) of Lemma \ref{modmulti}, the $g_i$ are allowed to
depend on $x$. When we apply this lemma in Theorem \ref{aboutissement}, we
do not make use of this. One might suspect that we could therefore
formulate Theorem \ref{aboutissement} in greater generality, namely that
we could permit the $b_i$ to vary locally. This seems to be a false
conclusion: The problem
seems to be that the identity in the hypotheses of Theorem
\ref{aboutissement} is an
identity in the ring $A$ whereas (i) in Lemma \ref{modmulti} is really on
the level of polynomials. If the $b_i$ depended on $x\in S(T)$, then also
the $g_i$ in equation (\ref{lacleman}) and we could not keep the left hand side
of (\ref{lacleman}) constant.
\end{remark}

\begin{remark}\label{passage}
One is tempted to think that, in the preceding proof, the passage from $T$
to the finitely generated semiring $T'\subset T$ would better have been carried
out already when choosing the finite empty subintersection of
(\ref{iplus})--(\ref{ia}). Though we could indeed have let run intersection
(\ref{it}) only over finitely many $t$ (analogously to intersections
(\ref{iplus}) and (\ref{itimes})), we then would not have known how to show
(\ref{malheure}) which was absolutely necessary to show $\ps\in S$.
\end{remark}

\begin{corollary}[Real Representation Theorem]\label{classique}
Let $T$ be a weakly divisible archimedean semiring of $A$. Suppose that
$a\in A$ satisfies $a>0$ on $S(T)$. Then $a\in T$.
\end{corollary}

\begin{proof}
Use $a=a\cdot 1$ as the required identity in the previous theorem.
\end{proof}

\begin{remark}\label{quark}
It is instructive to look how this section could be thinned out when one is
content with proving (rather than extending) the just stated Real
Representation Theorem.
The whole proof then collapses into what is essentially already contained in
the author's earlier work \cite{sw1} (see also \cite{sw3}). In the same way
than \cite{sw1} therefore can be read as a proof of the Real Representation
Theorem, the author's approach \cite[Section 2]{sw4} to
Putinar's Theorem \cite{put}\cite[Theorem 5.3.8]{pd} via P\'olya's theorem
(Corollary \ref{polya} above) can be read as a proof of Jacobi's variant
of the Real Representation Theorem \cite{jac}\cite[Theorem 5.3.6]{pd}.
Jacobi's variant says that Theorem \ref{classique} holds for
\emph{quadratic modules} au lieu of semirings where $T\subset A$ is
called a quadratic module if $0,1\in T$, $T+T\subset T$ and $A^2T\subset T$.
Scheiderer recently extended also this membership criterion of Jacobi from
positive to certain nonnegative elements \cite[Proposition 1.4]{s2} 
(see also \cite[p. 2, footnote 1]{m}). But the author's mentioned approach
via P\'olya's theorem to Jacobi's criterion seems not to be extendable
to this recent result of Scheiderer.
\end{remark}

\section{Alternative proof for preorders}\label{preordersection}

In this section, we demonstrate that Theorem \ref{aboutissement} can easily
be deduced from recent work of Scheiderer, Kuhlmann, Marshall and Schwartz
but only in the case where $T$ is a preorder. The following key lemma and its
proof is essentially \cite[Corollary 2.2]{kms}.

\begin{lemma}[Kuhlmann, Marshall, Schwartz]\label{primelemma}
Let $T$ be an archimedean preorder of $A$.
Suppose $1\in Aa+Ab$, $a,b\ge 0$ on $S(T)$ and $ab\in T$.
Then $a,b\in T$.
\end{lemma}

\begin{proof}
By our hypothesis and \cite[Lemma 2.1]{kms}
(see also \cite[Proposition 2.7]{s3} or \cite[Lemma 3.2]{m}
for a natural generalization of this not needed here), we have
$s,t\in A$ such that
$1=sa+tb$ and $s,t>0$ on $S(T)$. By the classical Real Representation Theorem
\ref{classique}, we have $s,t\in T$. Now $a=sa^2+tab\in T$ (here we use
that $A^2\subset T$). Symmetrically, we have of course $b\in T$.
\end{proof}

The next example shows that this key lemma does no longer hold in the general
situation where $T$ is only assumed to be a semiring instead of a preorder.

\begin{example}\label{exemplecontraire}
Let $A:=\R[X]$ and $T\subset A$ be the semiring generated by $\nn$ and the
three polynomials $1\pm X$ and $X^2+X^4$. The elements of $T$ are the
nonnegative linear combinations of products of these polynomials. By
Proposition \ref{archprop}, $T$ is clearly archimedean. Setting $a:=X^2$ and
$b:=1+X^2$, we clearly have $1\in Aa+Ab$ and $ab\in T$. Being sums of squares,
$a$ and $b$ are of course nonnegative on $S(T)$. We claim that $a\not\in T$.
Otherwise, we would have an identity
$$X^2=\sum_{\al\in\N^3}\la_\al(1+X)^{\al_1}(1-X)^{\al_2}(X^2+X^4)^{\al_3}
\qquad (\la_\al\ge 0).$$
Evaluating at $0$, we would get that the sum over all $\la_\al$
with $\al_3=0$ is $0$. But then, those $\la_\al$ would have to equal zero
since they are nonnegative. As a consequence, $X^2+X^4$ would divide $X^2$
which is absurd.
\end{example}

The idea for the next proof is from Corollaries 2.3 and 2.4 in \cite{kms}.

\begin{proof}[Alternative proof of Theorem \ref{aboutissement} in case
$A^2\subset T$]
The set $T':=T-a^2T\subset A$ is an archimedean preorder and we have
$S_{a=0}(T)=S(T')$. By hypothesis, we have therefore $b_i>0$ on $S(T')$ for
all $i$. From the classical Real Representation Theorem \ref{classique}, we
obtain $b_i\in T'$ for all $i$. Regarding the identity from the hypotheses,
this entails $a\in T'$, i.e., $a(1+at)\in T$
for some $t\in T$. By Lemma \ref{primelemma}, therefore $a\in T$.
\end{proof}

Even if Lemma \ref{primelemma} were true for semirings instead of preorders
(which is not the case), this alternative proof would break down. We would
have to replace the preordering $T'$ generated by $T$ and $-a^2$ by the
semiring $T-a^2T+a^4T-a^6T+\dots$ generated by $T$ and $-a^2$. But then we
would get only that
$$
a(1+at_1-a^3t_3+a^5t_5-a^7t_7+\dots)\in T\qquad\text{for some $t_1,t_3,\ldots
\in T$}
$$
instead of $a(1+at)\in T$ for some $t\in T$. The negative signs appearing in
the second factor of this product now prevent us from applying Lemma
\ref{primelemma}.

\section{Handelman's Theorem on powers of polynomials}\label{handelmansection}

In this section, we show that Theorem \ref{aboutissement} can be used to give
a new proof of a nice theorem of Handelman on powers of polynomials. See
Theorem \ref{handelman} and Corollary \ref{odd} below.
The original proof in \cite{h5} relies on some nontrivial facts from a whole
theory of a certain class of partially ordered abelian groups which is to a
large extent due to Handelman. Some of the used facts would not make sense in
our ring-theoretic setting, e.g., \cite[Proposition I.2(c)]{h3}. We have
decided to expose the whole material we need though a big part of it can be
found in less algebraic terminology in Handelman's original work
\cite{h1}\cite{h5} and in another new exposition of part of Handelman's
theory \cite{at}. This is not only because we want to keep this article
self-contained but also because we want to take on a new valuation theoretic
viewpoint. We will however only use the most basic facts and notions from
valuation theory as they can be found, for example, in the appendix of
\cite{pd}.

At first glance, it seems that our theorem is not suitable to prove Theorem
\ref{handelman}. Indeed, $\nnx$ is not an archimedean semiring of $\rx$.
However, for a semiring $T$ of a ring $A$, $T\cap O_T(A)$ is an archimedean
semiring of the ring of $T$-bounded elements $O_T(A)\subset A$ (cf. Lemma
\ref{archprop}). Still, this
does not seem to help since $O_{\nnx}(\rx)=\R$. When a ring of bounded elements
is too small, it is often a good idea to localize it by a fixed element, i.e.,
to build a new ring where division by this element is allowed
(see, e.g., \cite[Theorem 5.1]{sw2} or \cite{pv}). Following Handelman
(see, e.g., \cite[p. 61]{h3}), we will localize by a fixed
$0\neq g\in\nnx$. Hence we consider the ring
$$\rx_g:=\R\left[\x,\frac 1g\right]=
\left\{\frac f{g^k}\mid f\in\rx,k\in\N\right\}\subset\R(\x)$$
($\R(\x)$ denoting the quotient field of $\rx$) together with the
semiring
$$T_g:=\left\langle T\cup\left\{\frac 1g\right\}\right\rangle=
\left\{\frac f{g^k}\mid f\in\nnx,k\in\N\right\}\subset\rx_g$$
(we write angular brackets for the generated semiring).
For a polynomial $p\in\rx$, we denote by
$\text{Log}(p)\subseteq\N^n$ the set of all $\al\in\N^n$ for which the
coefficient of $\x^\al$ in $p$ does not vanish. Its convex hull
$\New(p)\subset\R^n$ is called the \emph{Newton polytope} of $p$.
It is easy to see that
\begin{align}
\label{logsub}
\Log(pq)&\subset\Log(p)+\Log(q)&&\text{for all $p,q\in\rx$,}\\
\label{logeq}
\Log(pq)&=\Log(p)+\Log(q)&&\text{for all $p,q\in\nnx$ and}\\
\label{neweq}
\New(pq)&=\New(p)+\New(q)&&\text{for all $p,q\in\rx$.}
\end{align}
These basic facts will frequently be used in the sequel, most often tacitly.
We now determine the ring of $T_g$-bounded elements $A(g)$ and its
(by Proposition \ref{archprop}) archimedean semiring $T(g):=T_g\cap A_g$:
\begin{align}\label{agdet}
A(g)&:=O_{T_g}(A_g)=\left\{\frac f{g^k}\mid f\in\rx,k\in\N,\Log(f)\subset
               \Log(g^k)\right\}\subset A_g\\
\label{tgdet}
T(g)&:=T_g\cap A(g)=\left\{\frac f{g^k}\mid f\in\nnx,k\in\N,\Log(f)\subset
               \Log(g^k)\right\}\subset A(g)
\end{align}
The inclusions from right to left are trivial whereas the inclusion
from left to right in (\ref{agdet}) uses (\ref{logsub}) and the one in
(\ref{tgdet}) uses (\ref{logsub}) and (\ref{logeq}). Using (\ref{logsub}),
the following becomes clear quickly:
\begin{align}\label{egendrer}
A(g)&=\R\left[\frac{\x^\al}g\mid\al\in\Log(g)\right]\\
T(g)&=\left\langle\nn\cup\left\{\frac{\x^\al}g\mid\al\in\Log(g)\right\}
      \right\rangle
\end{align}
Fix an arbitrary
$w\in\R^n$. There is exactly one valuation
$v_w:\R(\bar X)\to\R\cup\{\infty\}$ satisfying
\begin{equation}\label{vdef}
v_w(p)=-\max\{\langle w,\al\rangle\mid\al\in\text{Log}(p)\}
\qquad(0\neq p\in\rx).
\end{equation}
This is easy to show by noting that $\Log(p)$ can be replaced by $\New(p)$
in (\ref{vdef}) and using (\ref{neweq}).
Here and elsewhere $\langle w,\al\rangle$ denotes the usual scalar product
of $w$ and $\al$. We define the \emph{$w$-initial part}
$\init_w(p)\in\rx$ of a polynomial
$p\in\rx$ as the sum of those monomials appearing in $p$
belonging to an exponent tuple $\al\in\N^n$ for which $\langle w,\al\rangle$
gets maximal (i.e., equals $-v_w(p)$). The following is easy to check:
\begin{align}
\label{geoint}
\init_w(p)(x)&=\lim_{t\to\infty}e^{tv_w(p)}p(e^{tw_1}x_1,\dots,e^{tw_n}x_n)
&(0\neq p\in\rx,x\in\R^n)\\
\label{initmult}
\init_w(pq)&=\init_w(p)\init_w(q)&\text{($p,q\in\rx$)}
\end{align}
Let $\O_w$ denote the valuation ring  belonging to $v_w$ and
$\m_w$ its maximal ideal. It is an easy exercise to show that a ring
homomorphism
$\la_w:\O_w\to\R(\bar X)$ having kernel $\m_w$ is defined by
\begin{equation}\label{ladef}
\la_w\left(\frac pq\right):=\begin{cases}
                              0&\text{if\ } v_w(p)>v_w(q)\\
                              \frac{\text{in}_w(p)}{\text{in}_w(q)}
                              &\text{if\ } v_w(p)=v_w(q)
                              \end{cases}\qquad
                              (p,q\in\rx,q\neq 0),
\end{equation}
i.e., $\la_w$ is a place belonging to $v_w$.

We now give a concrete description of $S(T(g))$ using the notions just
defined. This result is from Handelman \cite[Theorem III.3]{h1} 
and also included in \cite[Lemma 2.4]{at}. For several reasons, we give here
a third exposition of this proof. In contrast to \cite[III.2]{h1} and
\cite[Lemma 2.3]{at}, we avoid the theory of polytopes and instead use some
basic valuation theory and (inspired by \cite[Lemma 1.10]{bra})
a fact from model theory. We believe that our viewpoint might be useful for
the investigation of rings other than $A(g)$.

\begin{theorem}[Handelman]\label{detspec}
For every $0\neq g\in\nnx$ and $x\in S(T(g))$, there is some $w\in\R^n$ and
$y\in\pos^n$ such that
$$a(x)=\la_w(a)(y)\qquad\text{for all $a\in A(g)$}.$$
\end{theorem}

\begin{proof}
By Chevalley's Theorem \cite[A.1.10]{pd}, we can extend the ring
homomorphism $x:A\to\R$ to a place of
$\R(\x)$, i.e., we find a valuation ring $\O\supset A(g)$ of $\R(\x)$ with
maximal ideal $\m$ and a ring homomorphism $\la:\O\to K$ into some extension
field $K$ of $\R$ with kernel $\m$ such that $\la|_{A(g)}=x$. Let
$v:\R(\x)\to\Ga\cup\{\infty\}$ be a valuation belonging to $\O$
where $\Ga$ is (after extension) without loss of generality a nontrivial
divisible ordered abelian group. Set
\begin{equation}\label{defla}
\La:=\left\{\al\in\Log(g)\mid v(\x^\al)=v(g)\right\}=
\left\{\al\in\Log(g)\mid\la\left(\frac{\x^\al}g\right)\neq 0\right\}.
\end{equation}
Now the first-order logic sentence
$$\exists u\exists v_1\dots\exists v_n\left(\bigwedge_{\al\in\La}
  \al_1v_1+\dots+\al_nv_n=u\wedge
  \bigwedge_{\al\in\Log(g)\setminus\La}
  \al_1v_1+\dots+\al_nv_n>u\right)$$
in the language $\{+,<,0\}$ holds in $\Ga$ (take $v(g)$ for $u$ and
$v(X_i)$ for $v_i$). It is a well-known fact in basic
model theory that all nontrivial divisible ordered abelian groups satisfy
exactly the same first-order sentences in this
language \cite[Corollary 3.1.17]{mar}. In particular, the above sentence
holds in $\R$, i.e., we find $w\in\R^n$ and $c\in\R$ such that
$\langle w,\al\rangle=c$ for all $\al\in\La$ and $\langle w,\al\rangle>c$
for all $\al\in\Log(g)\setminus\La$. It follows that $v_w(g)=-c$ and
\begin{equation}\label{rela}
\La=\{\al\in\Log(g)\mid v_w(\x^\al)=v_w(g)\}=
\left\{\al\in\Log(g)\mid\la_w\left(\frac{\x^\al}g\right)\neq 0\right\}.
\end{equation}
In view of (\ref{defla}), (\ref{rela}) and (\ref{egendrer}), it remains
only to show that there exists $y\in\pos^n$ such that
\begin{equation}\label{amontrer}
\la\left(\frac{\x^\al}g\right)=\la_w\left(\frac{\x^\al}g\right)(y)
\qquad\text{for all $\al\in\La$.}
\end{equation}
Now set $m:=\#\La -1\in\N$ and write $\La=\{\al^{(0)},\dots,\al^{(m)}\}$.
Assume for the moment that we have already shown the existence
of some $y\in\pos^n$ satisfying
\begin{equation}\label{amontrerii}
\la(\x^{\al^{(i)}-\al^{(0)}})=y^{\al^{(i)}-\al^{(0)}}
\qquad\text{for each $i\in\{1,\dots,m\}$.}
\end{equation}
Then we get immediately that even
\begin{equation}\label{iii}
\la(\x^{\al^{(i)}-\al^{(j)}})=y^{\al^{(i)}-\al^{(j)}}=
\la_w(\x^{\al^{(i)}-\al^{(j)}})(y)
\end{equation}
for $i,j\in\{0,\dots,m\}$.
Writing $g=\sum_{\al\in\Log(g)}a_\al\x^\al$, we obtain
\begin{multline*}
\la_w\left(\frac g{\x^{\al^{(i)}}}\right)(y)
\la\left(\frac{\x^{\al^{(i)}}}g\right)=
\sum_{\al\in\Log(g)}a_\al\la_w\left(\frac{\x^{\al}}{\x^{\al^{(i)}}}\right)(y)
\la\left(\frac{\x^{\al^{(i)}}}g\right)\\
\overset{(\ref{rela})}=\sum_{j=0}^ma_{\al^{(j)}}\la_w
\left(\frac{\x^{\al^{(j)}}}{\x^{\al^{(i)}}}\right)(y)
\la\left(\frac{\x^{\al^{(i)}}}g\right)
\overset{(\ref{iii})}=\sum_{j=0}^ma_{\al^{(j)}}\la
\left(\frac{\x^{\al^{(j)}}}{\x^{\al^{(i)}}}\right)
\la\left(\frac{\x^{\al^{(i)}}}g\right)\\
=\sum_{j=0}^ma_{\al^{(j)}}\la
\left(\frac{\x^{\al^{(j)}}}g\right)
\overset{(\ref{defla})}=\sum_{\al\in\Log(g)}a_\al
\la\left(\frac{\x^\al}g\right)=\la\left(\frac gg\right)=\la(1)=1
\end{multline*}
which shows (\ref{amontrer}). Therefore we are left with showing that
there is some $y\in\pos^n$ fulfilling (\ref{amontrerii}). Set
$\be^{(i)}:=\al^{(i)}-\al^{(0)}\in\Z^n$ and $z_i:=\la(\x^{\be^{(i)}})$ for
$i\in\{1,\dots,m\}$. Note that for all $i\in\{1,\dots,m\}$,
$$z_i=
\underbrace{\la\underbrace{\left(\frac{\x^{\al^{(i)}}}g\right)}_{\in T}}_
{\text{$\neq 0$ by (\ref{defla})}}
{\underbrace{\la\underbrace{\left(\frac{\x^{\al^{(0)}}}g\right)}_{\in T}}_
{\text{$\neq 0$ by (\ref{defla})}}}^{-1}>0
$$
since $\ph(T)\subseteq\nn$.
Using
$$y^{\be^{(i)}}=y_1^{\be_1^{(i)}}\dotsm y_n^{\be_n^{(i)}}=
  e^{(\log y_1)\be_1^{(i)}+\dots+(\log y_n)\be_n^{(i)}},$$
taking logarithms in (\ref{amontrerii}) and rewriting it in matrix form,
we therefore have to show that there are $y_1',\dots,y_n'\in\R$
(corresponding to $\log y_1,\dots,\log y_n$) such that
\begin{equation}\label{matrices}
  \underbrace{\begin{pmatrix}\log z_1&\dots&\log z_m\end{pmatrix}}_
  {=:L\in\R^{1\times m}}=
  \begin{pmatrix}y_1'&\dots&y_n'\end{pmatrix}
  \underbrace{\begin{pmatrix}\be_1^{(1)}&\dots&\be_1^{(m)}\\
                 \vdots&&\vdots\\
                 \be_n^{(1)}&\dots&\be_n^{(m)}
  \end{pmatrix}}_{=:B\in\R^{n\times m}}
\end{equation}
Provided now that $\ker B\subset\ker L$, the mapping
$\im B\to\R:Bv\mapsto Lv$ ($v\in\R^m$) is a well-defined linear map and can be
linearly extended to a map $\R^n\to\R$ represented by a $1\times n$
matrix $\begin{pmatrix}y_1'&\dots&y_n'\end{pmatrix}$
satisfying (\ref{matrices}).

Finally, we show $\ker B\subset\ker L$. Since all entries of $B$ lie in the
field $\Q$, $\ker B$ has a $\Q$-basis but then also $\R$-basis consisting of
vectors $k\in\Z^m$. Therefore consider an arbitrary $k\in\Z^m$ with
$$\sum_{j=1}^m\be_i^{(j)}k_j=0\qquad\text{for all $i\in\{1,\dots,m\}.$}$$
Taking the logarithm of
\begin{multline*}
e^{(\log z_1)k_1+\dots+(\log z_m)k_m}=z_1^{k_1}\dotsm z_m^{k_m}
=\la(\x^{\be^{(1)}})^{k_1}\dotsm\la(\x^{\be^{(m)}})^{k_m}\\
=\la(\x^{\be^{(1)}k_1+\dots+\be^{(m)}k_m})
=\la(\x^0)=\la(1)=1=e^0,
\end{multline*}
we get indeed $k\in\ker L$.
\end{proof}

\begin{corollary}[Handelman]\label{detspecii}
For every $0\neq g\in\nnx$ and $x\in S(T(g))$, there exist $w\in\R^n$ and
$y\in\pos^n$ such that
$$a(x)=\lim_{t\to\infty}a(e^{tw_1}y_1,\dots,e^{tw_n}y_n)\qquad
\text{for all $a\in A(g)$}.$$
\end{corollary}

\begin{proof}
Rewrite the last theorem using (\ref{geoint}) and (\ref{ladef}).
\end{proof}

We need a little number theoretic fact to make Proposition \ref{nombres}
below available.

\begin{lemma}
Suppose $l_1,l_2\in\N$ are relatively prime in $\Z$. Then for all $k\in\N$,
there exists
$m\in\N$ such that $\N\cap[m,\infty)\subset(\N\cap[k,\infty))l_1+
(\N\cap[k,\infty))l_2.$
\end{lemma}

\begin{proof}
Write $1=s_1l_1+s_2l_2$ with $s_1,s_2\in\Z$. If $s_1,s_2\ge 0$ then
either $s_1=l_1=1$ or $s_2=l_2=1$. Given $k\in\N$, we then can set $r:=k$.
Hence suppose, say, $s_1<0$. Then necessarily $l_2,s_2>0$. Given $k\in\N$, set
$$r:=(l_2-1)(-s_1)l_1+kl_1+kl_2\in\N.$$
Now we have for all $i\in\N$ and $j\in\{0,\dots,l_2-1\}$,
$$r+il_2+j=(k+(l_2-1-j)(-s_1))l_1+(k+js_2+i)l_2.$$
\end{proof}

\begin{proposition}\label{nombres}
Suppose $f\in\rx$ and let $l_1,l_2\in\N$ be relatively prime in $\Z$.
If it is true for $f^{l_1}$ and $f^{l_2}$ that all its sufficiently
high powers have nonnegative coefficients, then the same is true for $f$.
\end{proposition}

\begin{lemma}[Handelman]\label{gonflage}
Suppose $f\in\rx$, $1\le l\in\N$ and $f^l\in\rx^+$.
Then there is
$k_0\in\N$ such that for all $k\ge k_0$ and for all vertices $\al$
(i.e., extreme points) of $\New(f)$,
$$(lk-1)\al+\Log(f)\subset\Log(f^{lk}).$$
\end{lemma}

\begin{proof}
It is convenient to work in the ring
$\R[X_1,\dots,X_n,X_1^{-1},\dots,X_n^{-1}]\subset\R(\x)$ of Laurent
polynomials.
The Laurent monomials $\x^\al:=X_1^{\al_1}\dotsm X_n^{\al_n}$
($\al\in\Z^n$) form an $\R$-vector space
basis of it. Extending the definitions in the obvious way, we can speak
of $\Log(f)\subset\Z^n$ and $\New(f)\subset\R^n$ for any Laurent polynomial
$f$. We now prove our claim even for Laurent polynomials $f$.

Since the polytope $\New(f)$ has only finitely many vertices, it suffices
to show that the claimed inclusion of sets holds for a \emph{fixed}
vertex $\al$ and all large $k$.
Replacing $f$ by $\x^{-\al}f$, we can assume right away that $\al=0$.
Because the origin is now a vertex of $\New(f)$, we can choose $w\in\R^n$ such
that $\langle w,\be\rangle >0$ for all $0\neq\be\in\Log(f)$.
For all $0\neq\be,\ga,\de\in\Log(f)$ with $\be=\ga+\de$, in the equality
$\langle w,\be\rangle=\langle w,\ga\rangle+\langle w,\de\rangle$ the two
terms on the right hand side are then smaller than the left hand side.
We need the following consequence from this: Calling a nonzero element of
$\Log(f)$ an \emph{atom} if it is not a sum of two nonzero
elements of $\Log(f)$, every element of $\Log(f)$ can be written as a finite
sum of atoms (the origin being the sum of zero atoms). Since $\Log(f)$ is
finite, we can choose $k\in\N$ such that every element of $\Log(f)$ is a
sum of at most $k$ such atoms. On the other hand, because $f^l$ has
nonnegative coefficients, $\Log(f^{lk})$ consists of the sums of $k$ elements
of $\Log(f^l)$. Using $0\in\Log(f)$, it is enough to show that all atoms are
contained in $\Log(f^l)$. This is clear from the fact that an atom $\al$ can
can be written as a sum of $l$ elements from $\Log(f)$ only in a trivial way.
In fact, the coefficient of $\x^\al$ in $f^l$ is $l$ times the coefficient
of $\x^\al$ in $f$ and therefore nonzero.
\end{proof}

Now we are enough prepared to give a proof of Handelman's result based on
our membership criterion.

\begin{theorem}[Handelman]\label{handelman}
Let $f\in\rx$ be a polynomial such that $f^k$ has no
negative coefficients for some $k\ge 1$ and $f(1,1,\dots,1)>0$. Then
for all sufficiently large $k\in\N$, $f^k$ has no negative coefficients.
\end{theorem}

\begin{proof}
For any polynomial $p\in\rx$, we write $p^+$ for the sum of its monomials
with positive coefficients and $p^-$ for the negated sum of its monomials
with negative coefficients. So we always
have $p=p^+-p^-$, $p^+,p^-\in\rx^+$ and $\Log(p^+)\dot\cup\Log(p^-)=\Log(p)$.
First, we prove the theorem under the additional assumption
\begin{equation}\label{addass}
\init_w(f)\in\nnx\qquad
\text{for all $w\in\R^n$ with $\init_w(f)\neq f$.}
\end{equation}
By Lemma \ref{gonflage}, we can choose $k\in\N$ such that $g:=f^k$ has no
negative coefficients and
\begin{equation}\label{choix}
(k-1)\al+\Log(f)\subset\Log(g)\qquad
\text{for all vertices $\al$ of $\New(f)$.}
\end{equation}
Pick an arbitrary vertex $\al_0$ of
$\New(f)$. Then we have for all $N\in\N$,
\begin{equation}\label{atexample}
\begin{split}
a:=\frac{\x^{(k-1)\al_0}f}g&=
\underbrace{\left(1-N
\overbrace{\sum_\al\frac{\x^{(k-1)\al}f^+}g\cdot
                         \frac{\x^{(k-1)\al}f^-}g}^{=:c_1}\right)}_{=:b_1(N)}
\underbrace{\frac{\x^{(k-1)\al_0}f^+}g}_{=:t_1}\\
&+\underbrace{\left(N\overbrace{\sum_\al\frac{\x^{(k-1)\al}f^+}g\cdot
                       \frac{\x^{(k-1)\al}f^+}g}^{=:c_2}-1\right)}_{=:b_2(N)}
\underbrace{\frac{\x^{(k-1)\al_0}f^-}g}_{=:t_2}
\end{split}
\end{equation}
where the indices of summation run over all vertices $\al$ of $\New(f)$.
We will show that for $N$ sufficiently big, (\ref{atexample}) serves as an
identity like it is required in Theorem \ref{aboutissement}
which we are going to apply to the ring $A:=A(g)$ together with its archimedean
semiring $T:=T(g)$. To do this, first of all, observe that all fractions
appearing in (\ref{atexample}) lie in $A$ by (\ref{choix}).

Claim 1: $a>0$ on $\pos^n$. From $f^k=g\in\nnx$, it follows that $a^k>0$
on $\pos^n$. Using the continuity of $a$ on the connected space
$\pos^n$, we obtain either $a>0$ on $\pos^n$ or $a<0$ on $\pos^n$. The latter
can be excluded using the hypothesis $f(1,1,\dots,1)>0$

Claim 2: $a\ge 0$ on $S(T)$. This follows from Claim 1 and Corollary
\ref{detspecii}.

Claim 3: $c_1=0$ on $S_{a=0}(T)$. Let  $w\in\R^n$. According to Theorem
\ref{detspec}, we would have to show that $\la_w(a)(y)=0$ implies
$\la_w(c_1)(y)=0$ for all $y\in\pos^n$. In fact, we show that
$\la_w(c_1)\neq 0$ implies $\la_w(a)=a$ which is clearly more by Claim 1.
So suppose that $\la_w(c_1)\neq 0$. Then there is some vertex $\al$ of
$\New(f)$ with $v_w(\x^{(k-1)\al}f^+)=v_w(g)=v_w(\x^{(k-1)\al}f^-)$. This
implies $v_w(f^+)=v_w(f^-)$ whence $\init_w(f)\not\in\nnx$. From
(\ref{addass}), we now deduce $\init_w(f)=f$. This means that for all exponent
tuples $\be\in\N^n$ appearing in $f$, $\langle w,\be\rangle=-v_w(f)$ is
constant. Being
vertices of $\New(f)$, both $\al_0$ and $\al$ are among these $\be$. We obtain
therefore $v_w(\x^{(k-1)\al_0}f)=(k-1)v_w(\x^{\al_0})+v_w(f)=kv_w(f)
=v_w(f^k)=v_w(g)$. Exploiting the definition (\ref{ladef})
of $\la_w$ together with $\init_w(\x^{(k-1)\al_0}f)=\x^{(k-1)\al_0}\init_w(f)=
\x^{(k-1)\al_0}f$ and $\init_w(g)=\init_w(f^k)=\init_w(f)^k=f^k=g$, we see
that $\la_w(a)=a$.

Claim 4: $\New(f)=\New(f^+)$. Of course, we have
$\New(f)\supset\New(f^+)$
since $\Log(f)\supset\Log(f^+)$. For the other inclusion, it clearly suffices
to show that every vertex $\al$ of $\New(f)$, is contained in $\Log(f^+)$.
But for such a vertex $\al$, $\init_w(f)=\{\la\x^\al\}$ for some $\la\in\R^n$
and $w\in\R^n$. Except in the case where $f=\la\x^\al$, it follows from
(\ref{addass}) that $\la>0$ whence $\al\in\Log(f^+)$. If $f=\la\x^\al$,
then $\la>0$ follows from $f(1,1,\dots,1)>0$.

Claim 5: $c_2>0$ on $S(T)$. Let $w\in\R^n$. By
Theorem \ref{detspec}, $\la_w(c_2)(y)>0$ for all $y\in\pos^n$ is what we would
have to show.
By definition of $\la_w$ it is enough to show that
$\la_w(c_2)\neq 0$ since $\x^{(k-1)\al}f^+$ has no negative coefficients.
We obtain from Claim 4 that $v_w(f^+)=v_w(f)$.
Choose a vertex $\al$ of $\New(f)$ such that $v_w(f)=v_w(\x^\al)$.
Then $v_w(\x^{(k-1)\al}f^+)=(k-1)v_w(\x^\al)+v_w(f^+)=kv_w(f)=v_w(f^k)=v_w(g)$.
Therefore $\la_w(c_2)\neq 0$ as desired.

Regarded as a continuous real-valued function on the compact space $S(T)$,
$c_2$
is bounded from below by some positive real number by Claim 5. Consequently,
we can choose $N\in\N$ so large that $b_2(N)=Nc_2-1>0$ on the whole of $S(T)$,
in particular on $S_{a=0}(T)$. By Claim 3, we have that $b_1(N)=1-Nc_1=1>0$ on
$S_{a=0}(T)$. Of course, $t_1,t_2\in T$. Altogether, we can apply Theorem
\ref{aboutissement} and see that $a\in T$. By definition of $T=T(g)$, this
means that $g^m\x^{(k-1)\al_0}f\in\nnx$ for some $m\in\N$. Omitting
$\x^{(k-1)\al_0}$ does not change this fact, so that
$f^{km+1}=g^mf\in\nnx$. At the same time, of course, $f^{km}=g^m\in\nnx$.
Proposition \ref{nombres} yields now that all sufficiently high powers of
$f$ lie in $\nnx$.

Thus we have shown the theorem under the assumption (\ref{addass}). Now in the
general case, we proceed by induction on the number of monomials appearing in
$f$. The case where $f$ has only one monomial is trivial. Now suppose that
$f$ has at least two monomials. The hypothesis implies clearly that
\begin{equation}\label{krisperl}
f>0\qquad\text{on $\pos^n$}.
\end{equation}

Let $w\in\R^n$ such that $\init_w(f)$ has less monomials than $f$.
For some $k\ge 1$, $(\init_w(f))^k=\init_w(f^k)\in\rx^+$ by the hypotheses on
$f$. Evaluating this at $(1,1,\dots,1)$, we
see that $\init_w(f)$ does not vanish at this point. Moreover, it is
nonnegative at the same point by (\ref{geoint}) and
(\ref{krisperl}). Altogether, we can apply the induction hypothesis on
$\init_w(f)$ to get that $\init_w(f^k)=(\init_w(f))^k\in\rx^+$
for all large $k$.

Since $\{\init_w(f)\mid w\in\R^n\}$ is of course finite, this shows that we
find $k_0\in\N$ such that for any $k\ge k_0$ and $w\in\R^n$ with
$\init_w(f)\neq f$, $\init_w(f^k)\in\nnx$. This shows that
(\ref{addass}) is satisfied with
$f$ replaced by $f^k$ for any $k\ge k_0$
(note that $\init_w(f^k)\neq f^k$ implies trivially
$\init_w(f)\neq f$). In particular, we find $l_1,l_2\in\N$ that are relatively
prime in $\Z$ such that (\ref{addass}) holds with $f$ replaced by
$f^{l_1}$ and $f^{l_2}$, e.g., take $l_1:=k_0$ and $l_2:=k_0+1$.
By the special case of the theorem already proved, we
get that $(f^{l_1})^k$ and $(f^{l_2})^k$ have no negative coefficients for
all large $k$. According to Proposition \ref{nombres}, this means that all
sufficiently high powers of $f$ have only nonnegative coefficients.
\end{proof}

\begin{corollary}[Handelman]\label{odd}
If some odd power of a real polynomial in several variables has only
nonnegative coefficients, then so do all sufficiently high powers.
\end{corollary}

\end{document}